\documentclass[11pt]{amsart}

\usepackage{amssymb}
\usepackage{latexsym}
\usepackage{amsmath}

\begin{document}

\newtheorem{theorem}{Theorem}[section]
\newtheorem{proposition}{Proposition}[section]
\newtheorem{definition}{Definition}[section]
\newtheorem{corollary}{Corollary}[section]

\title[Quantum automorphism groups]{Quantum automorphism groups of vertex-transitive graphs of order $\leq$ 11}

\author{Teodor Banica}

\address{Departement of Mathematics, Universite Paul Sabatier, 118 route de Narbonne, 31062 Toulouse, France}

\email{banica@picard.ups-tlse.fr}

\author{Julien Bichon}

\address{Laboratoire de Mathematiques Appliquees, Universite de Pau et des Pays de l'Adour, IPRA, Avenue de l'universite, 64000 Pau, France}

\email{bichon@univ-pau.fr}



\begin{abstract}
We study quantum automorphism groups of vertex-transitive graphs having less than 11 vertices. With one possible exception, these can be obtained from cyclic groups ${\mathbb Z}_n$, symmetric groups $S_n$ and quantum symmetric groups $\mathcal Q_n$, by using various product operations.
The exceptional case is that of the Petersen graph,
and we present some questions about it.
\end{abstract}

\maketitle

\section*{Introduction}

A remarkable fact, discovered by Wang in \cite{wa}, is that the symmetric group $S_n$ has a quantum analogue $\mathcal Q_n$. For $n\geq 4$ this quantum group is bigger than $S_n$, and fits into Woronowicz's formalism in \cite{wo}.

The quantum group $\mathcal Q_n$ is best understood via its representation theory: with suitable definitions, it appears as the Tannakian realisation of the Temperley-Lieb algebra (\cite{ba2}). This elucidates a number of questions regarding the Cayley graph, fusion rules, amenability, etc. More generally, this puts $\mathcal Q_n$ into the framework of free quantum groups of Van Daele and Wang (\cite{vdw}), where a whole machinery, inspired by work of Gromov, Jones, Voiculescu, Wassermann, Weingarten is now available.

The study of $\mathcal Q_n$, and of free quantum groups in general, focuses now on more technical aspects: matrix models (\cite{bc}, \cite{bm}), ergodic actions (\cite{bdv}, \cite{va}), harmonic analysis (\cite{vv}, \cite{ve}).

The other thing to do is to study subgroups of $\mathcal Q_n$. This was started independently by the authors in \cite{ba1}, \cite{ba2} and \cite{bi1}, \cite{bi2}, and continued in the joint paper \cite{bb}. The notion that emerges from this work is that of quantum automorphism group of a vertex-transitive graph.

In this paper we describe quantum automorphism groups of vertex-transitive graphs having $n\leq 11$ vertices, with one graph omitted. This enhances previous classification work from \cite{ba1}, \cite{ba2}, \cite{bb}, where we have $n\leq 9$, also with one graph omitted.

Needless to say, in this classification project the value of $n$ is there only to show how far our techniques go.

The four main features of the present work are:

(1) Product operations. We have general decomposition results for Cartesian and lexicographic products. These are motivated by the graphs $\texttt{Pr}(C_5),\texttt{Pr}(K_5)$ and $C_{10}(4)$, which appear at $n=10$.

(2) The discrete torus. Here $n=9$. We prove that its quantum group is equal to its classical group, namely $S_3\wr{\mathbb Z}_2$. This answers a question left open in \cite{ba2}, \cite{bb}, and provides the first example of a graph having a usual wreath product as quantum symmetry group.

(3) Circulant graphs. It is known from \cite{ba1} that the $n$-cycle with $n\neq 4$ has quantum symmetry group $D_n$. This is extended in \cite{ba2} to a bigger class of circulant graphs. Here we further enlarge the list of such graphs, with an ad-hoc proof for $C_{10}(2)$, which appears at $n=10$.

(4) The Petersen graph. This appears at $n=10$, and the corresponding quantum group seems to be different from the known ones. 
Our other techniques do not apply here:
it cannot be written as a graph product, and is not a circulant graph.
Neither could we carry a direct analysis as in the torus case
because of the complexity 
of some computations. 
However we prove that the corresponding quantum group
is not isomorphic to $\mathcal Q_5$.
In other words, we might have here a ``new'' quantum group. However, we don't have a proof, and the question is left open.

As a conclusion, we have two questions:

(I) First is to decide whether the Petersen graph produces or not a new quantum group. If it does, this would probably change a bit the landscape: in the big table at the end, based on work since Wang's paper \cite{wa}, all quantum groups can be obtained from ${\mathbb Z}_n,S_n,\mathcal Q_n$.

(II) A good question is to try to characterize graphs having no quantum symmetry. This paper provides many new examples, and we have found some more by working on the subject, but so far we were unable to find a conceptual result here.

\smallskip

The paper is organized as follows. Sections 1, 2 are quite detailed preliminary sections, the whole paper, or at least the ideas involved in it, being intended to be accessible to non-specialists. 
Sections 3, 4, 5, 6 deal with different kinds of graphs, once again in a quite self-contained way. In Section 7 we present the classification result, in the form of a big, independent table. In the last section we present a technical result about the 
quantum group of the Petersen graph.

\section{Quantum permutation groups}

In this paper we use the following simplified version of Woronowicz's
compact quantum groups \cite{wo}, which is the only 
one we need when dealing with quantum symmetries of
classical finite spaces.

\begin{definition} A Hopf ${\mathbb C}^*$-algebra is  a ${\mathbb C}^*$-algebra $A$ with unit, endowed with morphisms
\begin{eqnarray*}
\Delta&:&A\to A\otimes A\cr
\varepsilon&:&A\to{\mathbb C}\cr
S&:&A\to A^{op}
\end{eqnarray*}
satisfying the usual axioms for a comultiplication, counit and antipode, along with the extra condition $S^2=id$.
\end{definition}

The more traditional terminology
for such an object is that of a 
"universal Hopf $\mathbb C^*$-algebra of Kac type''.
The universality condition refers to the fact that the counit
and antipode are assumed to be defined on the whole $\mathbb C^*$-algebra 
$A$ (in full generality, these are only defined on a 
dense Hopf $*$-subalgebra) and the Kac condition
refers to the condition $S^2= id$.

We warn the reader that the Hopf $\mathbb C^*$-algebras we consider here 
are not Hopf algebras in the usual sense
(the tensor product in the definition is a $\mathbb C^*$-tensor product). However, they possess canonically
defined dense Hopf $*$-subalgebras, from which they
can be reconstructed using the universal $C^*$-completion procedure.
See the survey paper \cite{mava}.

\medskip

The first example is with a compact group $G$. We can consider the algebra of continuous functions $A={\mathbb C}(G)$, with operations
\begin{eqnarray*}
\Delta(f)&=&(g,h)\to f(gh)\cr
\varepsilon(f)&=&f(1)\cr
S(f)&=&g\to f(g^{-1})
\end{eqnarray*}
where we use the canonical identification $A\otimes A={\mathbb C} (G\times G)$.

The second example is with a discrete group $\Gamma$. We have here the algebra $A={\mathbb C}^*(\Gamma)$, obtained from the usual group algebra ${\mathbb C} [\Gamma]$ by the universal $\mathbb C^*$-completion procedure, with operations
\begin{eqnarray*}
\Delta(g)&=&g\otimes g\cr
\varepsilon(g)&=&1\cr
S(g)&=&g^{-1}
\end{eqnarray*}
where we use the canonical embedding $\Gamma\subset A$.

In general, associated to an arbitrary Hopf ${\mathbb C}^*$-algebra $A$ are a compact quantum group $G$ and a discrete quantum group $\Gamma$, according to the following heuristic formulae:
$$A={\mathbb C}(G)={\mathbb C}^*(\Gamma)$$
$$G=\widehat{\Gamma}$$
$$\Gamma=\widehat{G}$$
These formulae are made into precise statements in the first section  of Woronowicz'
seminal paper \cite{wo}.
They are pieces of Pontryagin duality
for locally compact quantum groups, whose latest version is given in \cite{kuva}.

The compact quantum group morphisms are defined in the usual manner:
if $A= \mathbb C(G)$ and $B= \mathbb C(H)$ are Hopf $\mathbb C^*$-algebras,
a quantum group morphism $H \rightarrow G$ arises from a Hopf $\mathbb C^*$-algebra morphism $\mathbb C(G) \rightarrow \mathbb C(H)$, and we say that
$H$ is a quantum subgroup of $G$ if the corresponding 
morphism $\mathbb C(G) \rightarrow \mathbb C(H)$ is surjective.
We refer to \cite{wa0} for more details on the compact quantum 
group language.

\smallskip

A square matrix $u= (u_{ij}) \in M_n(A)$ is said to be multiplicative if 
$$
\Delta(u_{ij})=\sum u_{ik}\otimes u_{kj} \quad {\rm and} \quad
\varepsilon(u_{ij})=\delta_{ij}$$
Multiplicative matrices correspond to corepresentations
of the Hopf $\mathbb C^*$-algebra $A$, that is, to representations
of the compact quantum group $G$
with $A = \mathbb C(G)$.
Such a multiplicative matrix $u$ will also be interpreted
as a linear map  $\mathbb C^n \longrightarrow \mathbb C^n \otimes A$.

In this paper we are essentially interested in the
following special type of multiplicative matrices.  

\begin{definition}
A magic unitary matrix is a square matrix, 
all of whose entries are projections and all of whose rows and columns are partitions of unity.
\end{definition}

Here we say that a finite family of projections
is a partition of unity if these projections are pairwise orthogonal
and if their sum equals 1.

\smallskip

As a first example, consider a finite group $G$ acting on a finite set $X$. The characteristic functions
$$p_{ij}=\chi\{\sigma\in G\mid \sigma(j)=i\}$$
form a magic unitary matrix, because the corresponding sets form partitions of $G$, when $i$ or $j$ varies. We have the following formulae for ${\mathbb C}(G)$:

\begin{eqnarray*}
\Delta(p_{ij})&=&\sum p_{ik}\otimes p_{kj}\cr
\varepsilon(p_{ij})&=&\delta_{ij}\cr
S(p_{ij})&=&p_{ji}
\end{eqnarray*}
and therefore $p=(p_{ij})$ is a multiplicative matrix.

In the particular case of the symmetric group $S_n$ acting on $\{1,\ldots ,n\}$, the Stone-Weierstrass theorem shows that entries of $p$ generate ${\mathbb C}(S_n)$. This suggests the following construction, due to Wang (\cite{wa}).

\begin{definition} The $\mathbb C^*$-algebra
$A_s(n)$ is the universal ${\mathbb C}^*$-algebra generated by $n^2$ elements $u_{ij}$, with relations making $u$ into a magic unitary matrix, 
and with morphisms
\begin{eqnarray*}
\Delta(u_{ij})&=&\sum u_{ik}\otimes u_{kj}\cr
\varepsilon(u_{ij})&=&\delta_{ij}\cr
S(u_{ij})&=&u_{ji}
\end{eqnarray*}
as comultiplication, counit and antipode, making it into a 
Hopf $\mathbb C^*$-algebra. 
\end{definition}

This Hopf $\mathbb C^*$-algebra was discovered by Wang \cite{wa}.
The corresponding compact quantum group is denoted 
$\mathcal Q_n$ and we call it the quantum permutation
group or quantum symmetric group. 
This is motivated by the fact that 
the algebra $A_s(n)$ is the biggest Hopf $\mathbb C^*$-algebra coacting on 
the algebra $\mathbb C^n$, which is to say that the quantum group $\mathcal Q_n$
is the biggest one acting on 
$\{1,\ldots ,n\}$. The coaction $u : \mathbb C^n \longrightarrow \mathbb C^n \otimes A_s(n)$
is defined on Dirac masses by
$$u(\delta_i)=\sum \delta_j\otimes u_{ji}$$
and verification of axioms of coactions, as well as proof of universality, is by direct computation. See \cite{wa}.

We have a surjective morphism of Hopf ${\mathbb C}^*$-algebras
$$A_s(n)\to {\mathbb C}(S_n)$$
mapping $u_{ij}$ to $p_{ij}$ for any $i,j$. This morphism expresses the fact that the compact quantum group corresponding to $A_s(n)$ contains $S_n$.

This map is an isomorphism for $n=2,3$, as known from \cite{ba2}, \cite{wa}, and explained in section 3 below. At $n=4$ we have Wang's matrix
$$u=\begin{pmatrix}p&1-p&0&0\cr 1-p&p&0&0\cr 0&0&q&1-q\cr 0&0&1-q&q \end{pmatrix}$$
with $p,q$ free projections, which shows that 
there exists an epimorphism
$A_s(4) \to \mathbb C^*(\mathbb Z_2 * \mathbb Z_2)$
and hence
$A_s(n)$ is not commutative and is infinite dimensional. The same remains true for any $n\geq 4$.

\section{Quantum automorphism groups of graphs}

Consider a finite graph $X$. 
In this paper this means that we have a finite set of vertices, 
and certain pairs of distinct vertices are connected by unoriented edges.

It is convenient to assume that the vertex set is $\{1,\ldots ,n\}$.

\begin{definition}
The adjacency matrix of $X$ is the matrix
$$d\in M_n(0,1)$$
given by $d_{ij}=1$ if $i,j$ are connected by an edge, and $d_{ij}=0$ if not.
\end{definition}

The adjacency matrix is symmetric, and has $0$ on the diagonal. In fact, graphs having vertex set $\{1,\ldots ,n\}$ are in one-to-one correspondence with $n\times n$ symmetric 0--1 matrices having $0$ on the diagonal.

The quantum automorphism group of $X$ is obtained as an appropriate subgroup of the quantum permutation group of $\{1,\ldots ,n\}$. At level of Hopf $\mathbb C^*$-algebras, this means taking an appropriate quotient of $A_s(n)$.

\begin{definition}
Associated to a finite graph $X$ is the $\mathbb C^*$-algebra
$$A(X)=A_s(n)/<du=ud>$$
where $n$ is the number of vertices, and $d$ is the adjacency matrix.
\end{definition}

Since a permutation of the set $X$ is a graph automorphism
if and only if the corresponding permutation matrix
commutes with the adjacency matrix, it is reasonable
to say that  the quantum group corresponding to $A(X)$
is the quantum automorphism group of $X$.
In this way we have a commutative diagram of Hopf ${\mathbb C}^*$-algebras
$$\begin{matrix} 
A_s(n)&\ &\rightarrow&\ &A(X)\cr
 \ \cr \downarrow&\ &\ &\ 
  &\downarrow\cr \ \cr 
  {\mathbb C}(S_n)&\ &\rightarrow&\ &{\mathbb C}(G)
\end{matrix}$$
where $G=G(X)$ is the usual automorphism group of $X$, with the kernel of the right arrow being the commutator ideal of $A(X)$. Moreover, for a graph without edges we get indeed $A_s(n)$, and we have the formula
$$A(X)=A(X^c)$$
where $X^c$ is the complement of $X$. See \cite{ba2}, \cite{bb} for details.

The defining equations $ud = du$ of $A(X)$ means that
$d$, considered as a linear map $\mathbb C^n \rightarrow \mathbb C^n$,
is a morphism in the category of corepresentations of $A(X)$,
i.e. a morphism in the category of representations 
of the quantum group dual to $A(X)$. General properties of the
representation category of a compact quantum group
(see e.g. \cite{wo}) now ensure that the spectral projections
occurring in the spectral decomposition of $d$ are 
corepresentations morphisms, and hence the corresponding 
eigensubspaces are subcorepresentations. This key fact will be
used freely in the paper.

\smallskip

The following notion will play a central role in this paper.

\begin{definition}
We say that $X$ has no quantum symmetry if
$$A(X)={\mathbb C}(G)$$
where $G=G(X)$ is the usual automorphism group of $X$.
\end{definition}

This is the same as saying that $A(X)$ is commutative, because by the above considerations, ${\mathbb C}(G)$ is its biggest commutative quotient.

\medskip

We are particularly interested in the case of graphs $X$ having the property that $G$ acts transitively on the set of vertices.
These graphs were called homogeneous in previous work \cite{ba2}, \cite{bb},
but we use here the following more traditional terminology.

\begin{definition}
The graph $X$ is called vertex-transitive if for any two vertices $i,j$ there is $\sigma\in G(X)$  such that $\sigma(i)=j$.
\end{definition}

Each section of the paper ends with a small table, gathering information about vertex-transitive graphs having $\leq 11$ vertices. These small tables are to be put together in a single big table, at the end.

\smallskip

What we know so far is that we have
$$A(K_n)=A_s(n)$$
where $K_n$ is the complete graph having $n$ vertices. Moreover, we already mentioned that for $n=2,3$ the arrow
$$A_s(n)\to {\mathbb C}(S_n)$$
is an isomorphism, and for $n\geq 4$ it is not.

This information is summarized in the following table.

\begin{center}\begin{tabular}[t]{|l|l|l|l|l|}
\hline
Order & Graph & Classical group & Quantum group\\  
\hline \hline
2&$K_2$&$ {{\mathbb Z}_2}$&$ {{\mathbb Z}_2}$\\ 
\hline
3&$K_3$&$ {S_3}$&$ {S_3}$\\ 
\hline
$n\geq4$&$K_n$&$ S_n$&$ \mathcal Q_n$\\ 
\hline
\end{tabular}\end{center}

\medskip

Here in the right column $\mathcal Q_n$ with $n\geq 4$ is the compact quantum group associated to $A_s(n)$.

\section{Circulant graphs}

A graph with $n$ vertices is called circulant if its automorphism group contains a cycle of length $n$ (and hence in particular a 
copy of the cyclic group ${\mathbb Z}_n$). We are particularly interested in connected circulant graphs, which are the cycles with chords.

\begin{definition}
The graph $C_n(k_1,\ldots ,k_r)$, where
$$1<k_1<\ldots <k_r\leq [n/2]$$
are integers, is obtained by drawing the $n$-cycle $C_n$, then connecting all pairs of vertices at distance $k_i$, for any $i$.
\end{definition}

As basic examples, we have the $n$-cycle $C_n$, corresponding to the value $r=0$, and the $2n$-cycle with diagonals, $C_n^+=C_{2n}(n)$.

Observe that $K_n$ is a cycle with chords as well.

The adjacency matrix of a cycle with chords, denoted as usual $d$, is a circulant matrix. We use the following basic fact.

\begin{proposition}
We have $d(\xi^s)=2f(s)\xi^s$, where
$$f(s)=\sum_{i=0}^r \cos\left(\frac{2k_is\pi}{n}\right)
\quad \rm{(with} \  k_0=1\rm{)}$$
and $\xi$ is the vector whose coordinates are 
 $1,\omega , \ldots , \omega^{n-1}$ in the canonical basis of $\mathbb C^n$,
 with $\omega = e^{\frac{2i\pi}{n}}$.
\end{proposition}

This tells us that we have the following eigenspaces for $d$:
\begin{eqnarray*}
V_0&=&{\mathbb C} 1\cr
V_1&=&{\mathbb C}\xi +{\mathbb C}\xi^{n-1}\cr
V_2&=&{\mathbb C}\xi^2+{\mathbb C}\xi^{n-2}\cr
\ldots&&\ldots\cr
V_{m}&=&{\mathbb C}\xi^{m}+{\mathbb C}\xi^{n-m}\cr
\end{eqnarray*}
where $m=[n/2]$ and all sums are direct, except maybe for the last one, which depends on the parity of $n$.

The fact that these eigenspaces correspond or not to different eigenvalues depends of course on $f$.

We use the following result from \cite{ba2}, whose proof is briefly explained, because several versions of it will appear throughout the paper.

\begin{theorem}
If $n\neq 4$ and the associated function
$$f:\{1,2,\ldots,[n/2]\}\to{\mathbb R}$$
is injective, then $C_n(k_1,\ldots ,k_r)$ has no quantum symmetry.
\end{theorem}

\begin{proof}
Since ${\mathbb C}\xi\oplus{\mathbb C}\xi^{n-1}$ is invariant, the coaction can be written as
$${v}(\xi)=\xi\otimes a+\xi^{n-1}\otimes b$$
for some $a,b$. By taking powers and using $n\neq 4$ we get by induction
$${v}(\xi^s)=\xi^s\otimes a^s+\xi^{n-s}\otimes b^s$$
for any $s$, along with the relations $ab=-ba$ and $ab^2=ba^2=0$.

Now from ${v}(\xi^*)={v}(\xi)^*$ we get $b^*=b^{n-1}$, so $(ab)(ab)^*=ab^na^*=0$. Thus $ab=ba=0$, so $A(X)$ is commutative and we are done.
\end{proof}

\begin{corollary}
The following graphs have no quantum symmetry:
\begin{enumerate}
\item The cycles $C_n$ with $n\neq 4$.
\item The cycles with diagonals $C_8^+,C_{10}^+$.
\item The cycles with chords $C_9(3),C_{11}(2),C_{11}(3)$.
\end{enumerate}
\end{corollary}

\begin{proof} (1) follows from the fact that $f$ is decreasing, hence injective. As for (2) and (3), the corresponding 5 functions are given by
\begin{eqnarray*}
C_8^+&:&-0.29,1,-1.7,0\cr 
C_{10}^+&:&-0.19,1.3,-1.3,0.19,-2\cr
C_9(3)&:&0.26,-0.32,0.5,-1.43\cr
C_{11}(2)&:&1.25,-0.23,-1.10,-0.79,-0.11\cr
C_{11}(3)&:&0.69,-0.54,0.27,0.18.-1.61
\end{eqnarray*}
with $0.01$ error, so they are injective, and Theorem 3.1 applies.
\end{proof}

The graphs in Corollary 3.1 have usual symmetry group $D_n$, where $n$ is the number of vertices. We don't know if $G=D_n$ with $n\neq 4$ implies that the graph has no quantum symmetry. However, we are able to prove this for $n\leq 11$: graphs satisfying $G=D_n$ are those in Corollary 3.1, plus the cycle with chords $C_{10}(2)$, discussed below.

\begin{theorem}
The graph $C_{10}(2)$ has no quantum symmetry.
\end{theorem}

\begin{proof}
The function $f$ is given by
$$f(s)=\cos\left(\frac{s\pi}{5}\right)+\cos\left(\frac{2s\pi}{5}\right)$$
and we have $f(1)=-f(3)\simeq 1.11$, $f(2)=f(4)=-0.5$ and $f(5)=0$. Thus the list of eigenspaces is:
\begin{eqnarray*}
V_0&=&{\mathbb C} 1\cr
V_1&=&{\mathbb C}\xi\oplus{\mathbb C}\xi^{9}\cr
V_2&=&{\mathbb C}\xi^2\oplus{\mathbb C}\xi^4\oplus{\mathbb C}\xi^{6}\oplus{\mathbb C}\xi^{8}\cr
V_3&=&{\mathbb C}\xi^3\oplus{\mathbb C}\xi^{7}\cr
V_5&=&{\mathbb C}\xi^5
\end{eqnarray*}

Since coactions preserve eigenspaces, we can write
$${v}(\xi)=\xi\otimes a+\xi^{9}\otimes b$$
for some $a,b$. Taking the square of this relation gives
$${v}(\xi^2)=\xi^2\otimes a^2+\xi^{8}\otimes b^2+1\otimes (ab+ba)$$
and once again since ${v}$ preserves eigenspaces, we get $ab=-ba$. Taking now the cube of the above relation gives
\begin{eqnarray*}
{v}(\xi^3)&=&\xi^3\otimes a^3+\xi^{7}\otimes b^3
+\xi\otimes ba^2 +\xi^{9}\otimes ab^{2}
\end{eqnarray*}
and once again since ${v}$ preserves eigenspaces, we get:
$$ab^2=0= ba^2$$
With the relations $ab=-ba$ and $ab^2=ba^2=0$ in hand, we get by induction the formula
$${v}(\xi^s)=\xi^s\otimes a^s+\xi^{n-s}\otimes b^s$$
and we can conclude by using adjoints, as in proof of Theorem 3.1.
\end{proof}

For graphs having $n\leq 11$ vertices, results in this section are summarized in the following table.

\begin{center}\begin{tabular}[t]{|l|l|l|l|l|}
\hline
Order & Graph & Classical group & Quantum group\\  
\hline \hline
$n\geq5$&$C_n$&$ D_n$&$ D_n$\\ 
\hline
8&$C_8,C_8^+$&$ D_8$&$ D_8$\\ 
\hline
9&$C_9,C_9(3)$&$ D_9$&$ D_9$\\ 
\hline
10&$C_{10},C_{10}(2),C_{10}^+$&$ D_{10}$&$ D_{10}$\\ 
\hline
11&$C_{11},C_{11}(2),C_{11}(3)$&$ D_{11}$&$ D_{11}$\\ 
\hline
\end{tabular}\end{center}

\medskip

As already mentioned, we don't know if these computations are particular cases of some general result.

\section{Products of graphs}

For a finite graph $X$, it is convenient to use the notation
$$X=(X,\sim)$$
where the $X$ on the right is the set of vertices, and where we write $i\sim j$ when two vertices $i,j$ are connected by an edge.

\begin{definition}
Let $X,Y$ be two finite graphs.
\begin{enumerate}
\item 
The direct product $X\times Y$ has vertex set $X\times Y$, and edges
$$(i,\alpha)\sim(j,\beta)\Longleftrightarrow i\sim j,\, \alpha\sim\beta.$$
\item
The Cartesian product $X\square Y$ has vertex set $X\times Y$, and edges
$$(i,\alpha)\sim(j,\beta)\Longleftrightarrow i=j,\, \alpha\sim\beta\mbox{ \rm{or} }i\sim j,\alpha=\beta.$$
\end{enumerate}
\end{definition}

The direct product is the usual one in a categorical sense. As for the Cartesian product, this is a natural one from a geometric viewpoint:  for instance a product by a segment gives a prism.

\begin{definition}
The prism having basis $X$ is $\mathtt{Pr}(X)=K_2\square X$.
\end{definition}

We have embeddings of usual symmetry groups
$$G(X) \times G(Y) \subset G(X \times Y)$$
$$G(X) \times G(Y) \subset G(X \square Y)$$
which have the following quantum analogues.

\begin{proposition}
We have surjective morphisms of Hopf ${\mathbb C}^*$-algebras
$$A(X \times Y) \longrightarrow A(X) \otimes A(Y)$$
$$A(X \square Y) \longrightarrow A(X) \otimes A(Y).$$
\end{proposition}

\begin{proof}
We use the canonical identification 
$${\mathbb C}(X \times Y)={\mathbb C}(X) \otimes {\mathbb C}(Y)$$
given by $\delta_{(i,\alpha)}=\delta_i\otimes\delta_\alpha$. The adjacency matrices are given by
$$d_{X \times Y} = d_X \otimes d_Y$$
$$d_{X \square Y} = d_X \otimes 1 + 1 \otimes  d_Y$$
so if $u$ commutes with $d_X$ and $v$ commutes with $d_Y$, the matrix
$$u\otimes v=(u_{ij}v_{\alpha\beta})_{(i\alpha,j\beta)}$$
is a magic unitary that 
commutes with both $d_{X\times Y}$ and $d_{X\square Y}$. This gives morphisms as in the statement, and surjectivity follows by summing over $i$ and $\beta$.
\end{proof}

\begin{theorem}
Let $X$ and $Y$ be finite connected regular graphs. If their spectra $\{\lambda\}$ and $\{\mu\}$ do not contain $0$ and satisfy
$$\{ \lambda_i/\lambda_j\} \cap \{\mu_k/\mu_l\}
= \{1\}$$
then $A(X \times Y)=A(X) \otimes A(Y)$. Also, if their spectra satisfy
$$\{\lambda_i - \lambda_j \} \cap \{\mu_k - \mu_l\}
= \{0\}$$
then $A(X \square Y)=A(X) \otimes A(Y)$.
\end{theorem}

\begin{proof}
We follow \cite{ba2}, where the first statement is proved. Let $\lambda_1$ be the valence of $X$. Since $X$ is regular
we have $\lambda_1 \in {\rm Sp}(X)$, with $1$ as eigenvector,
and since $X$ is connected $\lambda_1$ has multiplicity one.
Hence if $P_1$ is the orthogonal projection onto 
${\mathbb C}1$, the spectral decomposition of $d_X$ is of the following form:
$$d_X = \lambda_1 P_1 + \sum_{i\not=1}\lambda_i P_i$$
We have a similar formula for $d_Y$:
$$d_Y = \mu_1 Q_1 + \sum_{j\not=1}\mu_j Q_j$$
This gives the following formulae for products:
$$d_{X\times Y}=\sum_{ij}(\lambda_i\mu_j)P_{i}\otimes Q_{j}$$
$$d_{X \square Y} = \sum_{i,j}(\lambda_i + \mu_i)P_i \otimes Q_j$$
Here projections form partitions of unity, and the scalar are distinct, so these are spectral decomposition formulae. We can conclude as in \cite{ba2}. 
The universal coactions will commute with any of the spectral
projections, and hence with both $P_1 \otimes 1$ and $1 \otimes Q_1$.   
In both cases the universal coaction $v$ is the tensor product of 
its restrictions to the images of $P_1\otimes 1$ 
(i.e. $1 \otimes \mathbb C(Y)$) and of $1\otimes Q_1$ 
(i.e. $\mathbb C(X) \otimes 1$).
\end{proof}

\begin{corollary}
\ 
\begin{enumerate}
\item We have $A(K_m \times K_n)=A(K_m) \otimes A(K_n)$ for $m \not = n$.
\item We have $A(\mathtt{Pr}(K_n))={\mathbb C}({{\mathbb Z}_2})\otimes A_s(n)$, for any $n$.
\item We have $A(\mathtt{Pr}(C_n))={\mathbb C}(D_{2n})$, for $n$ odd.
\item We have $A(\mathtt{Pr}(C_4))={\mathbb C}({{\mathbb Z}_2})\otimes A_s(4)$.
\end{enumerate}
\end{corollary}

\begin{proof}
The spectra of graphs involved are ${\rm Sp}(K_2)=\{-1,1\}$ and
\begin{eqnarray*}
{\rm Sp}(K_n)&=&\{ -1,\ n-1\}\cr
{\rm Sp}(C_n)&=&\{2\cos (2k\pi /n)\mid k=1,\ldots ,n\}
\end{eqnarray*}
so the first three assertions follow from Theorem 4.1. We have
$$\texttt{Pr}(C_4)=K_2\times K_4$$
(this graph is the cube)
and the fourth assertion follows from the first one.
\end{proof}

We get the following table, the product operation
$\times$ on quantum groups being the one dual
to the tensor product of Hopf $\mathbb C^*$-algebras.

\begin{center}\begin{tabular}[t]{|l|l|l|l|l|}
\hline
Order & Graph & Classical group & Quantum group\\  
\hline \hline
8&$\texttt{Pr}(C_4)$&$S_4\times{{\mathbb Z}_2}$&$\mathcal Q_4\times{{\mathbb Z}_2}$\\
\hline
10&$\texttt{Pr}(C_5)$&$D_{10}$&$D_{10}$\\
\hline
10&$\texttt{Pr}(K_5)$&$S_5\times{{\mathbb Z}_2}$&$\mathcal Q_5\times{{\mathbb Z}_2}$\\
\hline
\end{tabular}\end{center}

\section{The torus graph}

Theorem 4.1 doesn't apply to the case $X=Y$, and the problem of computing the algebras $A(X\times X)$ and $A(X\square X)$ appears.
At level of classical symmetry groups, there is no simple formula describing $G(X\times X)$ and $G(X\square X)$. Thus we have reasons to believe that the above problem doesn't have a simple solution either.

A simpler question is to characterize graphs $X$ such that $X\times X$ or $X\square X$ has no quantum symmetry. We don't have a general result here, but we are able however to deal with the case $X=K_3$.

\begin{definition}
The graph $\mathtt{Torus}$ is the graph $K_3\times K_3=K_3\square K_3$.
\end{definition}

The result below answers a question asked in \cite{ba2}, \cite{bb}. It also provides the first example of graph having a classical wreath product as quantum symmetry group.

\begin{theorem}
The graph $\mathtt{Torus}$ has no quantum symmetry. 
\end{theorem}

\begin{proof}
The spectrum of $K_3$ is known to be
$${\rm Sp}(K_3)=\{ -1,2\}$$
with corresponding eigenspaces given by
\begin{eqnarray*}
F_2&=&{\mathbb C} 1\cr
F_{-1}&=&{\mathbb C}\xi\oplus{\mathbb C}\xi^2
\end{eqnarray*}
where $\xi$ is the vector formed by third roots of unity.

Tensoring the adjacency matrix of $K_3$ with itself gives
$${\rm Sp}(\texttt{Torus}) = \{-2,1,4\}$$
with corresponding eigenspaces given by
\begin{eqnarray*}  
E_4&=&{\mathbb C}\xi_{00}\cr
E_{-2}&=&{\mathbb C}\xi_{10}\oplus {\mathbb C}\xi_{01}\oplus {\mathbb C}\xi_{20} \oplus {\mathbb C}\xi_{02}\cr
E_{1}&=&{\mathbb C}\xi_{11} \oplus{\mathbb C}\xi_{12} \oplus{\mathbb C}\xi_{21} \oplus{\mathbb C}\xi_{22}
\end{eqnarray*}
where we use the notation $\xi_{ij}=\xi^i\otimes \xi^j$.

The universal coaction $v$ preserves eigenspaces, so we have
\begin{eqnarray*}
v(\xi_{10})&=& \xi_{10} \otimes a + \xi_{01} \otimes b + \xi_{20} \otimes c
+ \xi_{02} \otimes d\cr
v(\xi_{01})&=&\xi_{10} \otimes\alpha + \xi_{01} \otimes \beta +
\xi_{20} \otimes \gamma + \xi_{02} \otimes \delta
\end{eqnarray*}
for some $a,b,c,d,\alpha,\beta,\gamma,\delta \in A$. Taking the square of $v(\xi_{10})$ gives 
$$v(\xi_{20})=\xi_{20} \otimes a^2 + \xi_{02} \otimes b^2 + \xi_{10} \otimes c^2
+\xi_{01} \otimes d^2$$
along with relations coming from eigenspace preservation:
$$ab = -ba, \ ad=-da, \ bc = -cb , \ cd = -dc$$
$$ac+ca = -(bd+db)$$

Now since $a,b$ anticommute, their squares have to commute.
On the other hand, by applying $v$ to the equality $\xi_{10}^*=\xi_{20}$, we get the following formulae for adjoints:
$$a^* = a^2, \ b^*=b^2, \ c^* = c^2, \ d^* = d^2$$

The commutation relation $a^2b^2=b^2a^2$ reads now $a^*b^*=b^*a^*$, and by taking adjoints we get $ba=ab$. Together with $ab=-ba$ this gives:
$$ab=ba=0$$

The same method applies to $ad,bc,cd$, and we end up with:
$$ab=ba=0,\ ad=da =0, \ bc=cb =0, \ cd = dc =0$$

We apply now $v$ to the equality $1=\xi_{10}\xi_{20}$. We get that $1$ is the sum of $16$ terms, all of them of the form $\xi_{ij}\otimes P$, where $P$ are products between $a,b,c,d$ and their squares. Due to the above formulae 8 terms vanish, and the $8$ remaining ones produce the formula
$$1=a^3 +b^3 +c^3 +d^3$$
along with relations coming from eigenspace preservation:
$$ac^2=ca^2=bd^2=db^2=0$$

Now from $ac^2=0$ we get $a^2c^2=0$, and by taking adjoints this gives $ca=0$. The same method applies to $ac,bd,db$, and we end up with:
$$ac=ca=0,\ bd=db=0$$

In the same way one shows that $\alpha,\beta,\gamma,\delta$
pairwise commute:
$$\alpha\beta=\beta\alpha=\ldots =\gamma\delta=\delta\gamma=0$$

It remains to show that $a,b,c,d$ commute with $\alpha,\beta,\gamma,\delta$. For, we apply $v$ to the following equality:
$$\xi_{10}\xi_{01}=\xi_{01}\xi_{10}$$

We get an equality between two sums having 16 terms each, and by using 
again eigenspace preservation we get the following formulae relating the corresponding 32 products $a\alpha, \alpha a$ etc.:
$$a\alpha = 0 = \alpha a , \ b\beta =0 = \beta b , \
c\gamma = 0 = \gamma c , \ d \delta =0 = \delta d,$$
$$a\gamma + c\alpha + b \delta + d\beta = 0 = 
\alpha c + \gamma a + \beta d + \delta b,$$
$$a\beta +b \alpha = \alpha b + \beta a, \
b \gamma + c\beta =  \beta c + \gamma b ,$$ 
$$c\delta + d \gamma =  \gamma d + \delta c , \
a \delta + d\alpha = \alpha d + \delta a$$  

Multiplying the first equality in the second row on the left by $a$ and on the right by $\gamma$ gives $a^2\gamma^2 =0$, and by taking adjoints we get $\gamma a=0$. The same method applies to the other 7 products involved in the second row, so all 8 products involved in the second row vanish:
$$a\gamma =c\alpha=b\delta = d\beta= \alpha c=\gamma a
=\beta d=\delta b=0$$

We use now the first equality in the third row. Multiplying it on the left by $a$ gives $a^2\beta=a\beta a$, and multiplying it on the right by $a$ gives $a\beta a=\beta a^2$. Thus we get the commutation relation $a^2\beta=\beta a^2$.

On the other hand from $a^3+b^3+c^3+d^3=1$ we get $a^4=a$, so:
$$a\beta = a^4 \beta = a^2a^2 \beta = \beta a^2 a^2 = \beta a$$

One shows in a similar manner that the missing commutation formulae $a\delta = \delta a$ etc. hold as well. Thus $A$ is commutative.    
\end{proof}

\begin{center}\begin{tabular}[t]{|l|l|l|l|l|}
\hline
Order & Graph & Classical group & Quantum group\\  
\hline \hline
9&$\texttt{Torus}$&$S_3\wr{{\mathbb Z}_2}$&$S_3\wr{{\mathbb Z}_2}$\\
\hline
\end{tabular}\end{center}

\medskip

\section{Lexicographic products}

Let $X$ and $Y$ be two finite graphs. Their lexicographic product is obtained by putting a copy of $X$ at each vertex of $Y$:

\begin{definition}
The lexicographic product $X\circ Y$ has vertex set $X\times Y$, and edges are given by
$$(i,\alpha)\sim(j,\beta)\Longleftrightarrow \alpha\sim\beta\mbox{ \rm{or} }\alpha=\beta,\,
 i\sim j.$$
\end{definition}

The terminology comes from a certain similarity with the ordering of usual words, which is transparent when iterating $\circ$.

The simplest example 
is with  $X\circ X_n$, where $X_n$
is the graph having $n$ vertices and no edges: 
the graph $X\circ X_n$ is the graph consisting of $n$ disjoint copies of $X$.

\begin{definition}
$nX$ is the disjoint union of $n$ copies of $X$.
\end{definition}

When $X$ is connected, we have an isomorphism
$$G(nX)=G(X)\wr S_n$$
where $\wr$ is a wreath product. In other words, we have:
$$G(X\circ X_n)=G(X)\wr G(X_n)$$

In the general case, we have the following embedding of usual symmetry groups:
$$G(X)\wr G(Y)\subset G(X\circ Y)$$

The quantum analogues of these results use the notion of free wreath product from \cite{bi2, bb}. In the following definition, a pair
$(A,u)$ is what we call a quantum permutation group in \cite{bb}:
$A$ is a Hopf $\mathbb C^*$-algebra 
and $u$ is a multiplicative magic unitary matrix.

\begin{definition}
The free wreath product of $(A,u)$ and $(B,v)$ is
$$A*_wB=(A^{*n}*B)/<[u_{ij}^{(a)},v_{ab}]=0>$$
where $n$ is the size of $v$, with magic unitary matrix $w_{ia,jb}=u_{ij}^{(a)}v_{ab}$.
\end{definition}

In other words, $A*_wB$ is the universal $\mathbb C^*$-algebra generated by $n$ copies of $A$ and a copy of $B$, with the $a$-th copy of $A$ commuting with the $a$-th row of $v$, for any $a$. The Hopf $\mathbb C^*$-algebra structure on $A *_w B$
is the unique one making $w$ into a multiplicative matrix.  
With this definition, we have the following result (\cite{bb}).

\begin{theorem}
If $X$ is connected we have $A(nX)=A(X)*_wA_s(n)$.
\end{theorem}

Note that the embedding $A(X)^{*n}\hookrightarrow A(X)*_wA_s(n)$ ensures
that $A(nX)$ is an infinite-dimensional algebra whenever 
$n \geq 2$ and $G(X)$ is non trivial.

In the general case, 
we have the following quantum analogue of the embedding result for $G(X)\wr G(Y)$.

\begin{proposition}
We have a surjective morphism of Hopf ${\mathbb C}^*$-algebras
$$A(X\circ  Y) \longrightarrow A(X) *_w A(Y).$$
\end{proposition}

\begin{proof}
We use the canonical identification
$${\mathbb C}(X \times Y)={\mathbb C}(X) \otimes {\mathbb C}(Y)$$
given by $\delta_{(i,\alpha)}=\delta_i\otimes\delta_\alpha$. The adjacency matrix of $X\circ Y$ is
$$d_{X\circ Y} = d_X \otimes 1 + \mathbb I \otimes d_Y$$
where $\mathbb I$ is the square matrix filled with $1$'s.

Let $u,v$ be the magic unitary matrices of $A(X),A(Y)$. The magic unitary matrix of $A(X)*_wA(Y)$ is given by
$$w_{ia,jb}= u_{ij}^{(a)}v_{ab}$$
and from the fact that $u$ commutes with $d_X$ (and $\mathbb I$)
and $v$ commutes with $d_Y$, we get that $w$ commutes with $d_{X\circ Y}$. This gives a morphism as in the statement, and surjectivity follows by summing over $i$ and $b$. 
\end{proof}

\begin{theorem}
Let $X,Y$ be regular graphs, with $X$ connected. If their spectra $\{\lambda_i\}$ and $\{\mu_j\}$ satisfy the condition
$$\{ \lambda_1-\lambda_i\mid i\neq 1 \} \cap \{-n\mu_j\} = \emptyset$$
where $n$ and $\lambda_1$ are the order and valence of $X$, then $A(X \circ  Y)=A(X) *_w A(Y)$.   
\end{theorem}

\begin{proof}
We denote by $P_i,Q_j$ the spectral projections corresponding to $\lambda_i,\mu_j$. Since $X$ is connected we have $P_1=\frac{1}{n}\,{\mathbb I}$, and we get:
\begin{eqnarray*}
d_{X\circ Y}
&=&d_X\otimes 1+{\mathbb I}\otimes d_Y\cr
&=&\left(\sum_i\lambda_iP_i\right)\otimes\left(\sum_jQ_j\right)+\left(nP_1\right)\otimes \left(\sum_i\mu_jQ_j\right)\cr
&=&\sum_j(\lambda_1+n\mu_j)(P_1 \otimes Q_j) + \sum_{i\not=1}\lambda_i (P_i\otimes 1)
\end{eqnarray*} 

In this formula projections form a partition of unity and scalars are distinct, so this is the spectral decomposition of $d_{X\circ Y}$.

Let $W$ be the universal coaction on
$X\circ Y$. Then $W$ must commute with all spectral projections, and in particular:
$$[W,P_1 \otimes Q_j]=0$$

Summing over $j$ gives $[W, P_1 \otimes 1]=0$, so $1\otimes {\mathbb C}(Y)$ is invariant under the coaction. The corresponding restriction of $W$ gives a coaction of $A(X\circ Y)$
on $1 \otimes {\mathbb C}(Y)$, say
$$W(1 \otimes e_a) = \sum_b 1 \otimes e_b \otimes y_{ba}$$
where $y$ is a magic unitary. On the other hand we can write
$$W(e_i \otimes 1) = \sum_{jb} e_j \otimes e_b \otimes x_{ji}^b$$  
and by multiplying by the previous relation we get:

\begin{eqnarray*}
W(e_i \otimes e_a)
&=&\sum_{jb} e_j \otimes e_b \otimes
y_{ba}x_{ji}^b\cr
&=& \sum_{jb} e_j \otimes e_b \otimes x_{ji}^b y_{ba}
\end{eqnarray*}

This shows that coefficients of $W$ have the following form:
$$W_{jb,ia} = y_{ba} x_{ji}^b=x_{ji}^b y_{ba}$$

Consider now the matrix $x^b=(x_{ij}^b)$. Since $W$ is a morphism of algebras, each row of $x^b$ is a partition  of unity. Also using the antipode, we have
\begin{eqnarray*}
S\left(\sum_jx_{ji}^{b}\right)
&=&S\left(\sum_{ja}x_{ji}^{b}y_{ba}\right)\cr
&=&S\left(\sum_{ja}W_{jb,ia}\right)\cr
&=&\sum_{ja}W_{ia,jb}\cr
&=&\sum_{ja}x_{ij}^ay_{ab}\cr
&=&\sum_ay_{ab}\cr
&=&1
\end{eqnarray*}
so we conclude that $x^b$ is a magic unitary.

We check now that $x^a,y$ commute with $d_X,d_Y$. We have
$$(d_{X\circ Y})_{ia,jb} = (d_X)_{ij}\delta_{ab} + (d_Y)_{ab}$$
so the two products between $W$ and $d_{X\circ Y}$ are given by:
$$(Wd_{X\circ Y})_{ia,kc}=\sum_j W_{ia,jc} (d_X)_{jk} + \sum_{jb}W_{ia,jb}(d_Y)_{bc}$$
$$(d_{X\circ Y}W)_{ia,kc}=\sum_j (d_X)_{ij} W_{ja,kc} + \sum_{jb}(d_Y)_{ab}W_{jb,kc}$$

Now since $W$ commutes with $d_{X\circ Y}$, the terms on the right are equal, and by summing over $c$ we get:

$$\sum_j x_{ij}^a(d_X)_{jk} + \sum_{cb} y_{ab}(d_Y)_{bc}
= \sum_{j} (d_X)_{ij}x_{jk}^a + \sum_{cb} (d_Y)_{ab}y_{bc}$$

The graph $Y$ being regular, the second sums in both terms are equal to the valency of $Y$, so we get $[x^a,d_X]=0$.

Now once again from the formula coming from commutation of $W$ with $d_{X\circ Y}$, we get $[y,d_Y] =0$.

Summing up, the coefficients of $W$ are of the form
$$W_{jb,ia}=x_{ji}^by_{ba}$$
where $x^b$ are magic unitaries commuting with $d_X$, and $y$ is a magic unitary commuting with $d_Y$. This gives a morphism
$$A(X)*_wA(Y) \longrightarrow A(X\circ Y)$$
mapping $u_{ji}^{(b)}\to x_{ji}^b$
and $v_{ba}\to y_{ba}$, which is inverse to the morphism in the previous proposition.
\end{proof}

\begin{corollary} 
We have $A(C_{10}(4))= {\mathbb C}({\mathbb Z}_2)*_w{\mathbb C}(D_5)$.  
\end{corollary}

\begin{proof}
We have isomorphisms
$$C_{10}(4)=C_{10}(4,5)^c=K_2\circ C_5$$
and Theorem 6.2 applies to the product on the right.
\end{proof}

Together with Theorem 6.1, this corollary gives the following table,
where ${\,\wr_*\,}$ is defined by ${\mathbb C}(G{\,\wr_*\,}H)={\mathbb C}(G)*_w{\mathbb C}(H)$.

\begin{center}\begin{tabular}[t]{|l|l|l|l|l|}
\hline
Order & Graph & Classical group & Quantum group\\  
\hline \hline
4&$2K_2$&${\mathbb Z}_2\wr{\mathbb Z}_2$&${{\mathbb Z}_2}{\,\wr_*\,}{{\mathbb Z}_2}$\\
\hline
6&$2K_3$&${S_3}\wr{{\mathbb Z}_2}$&${S_3}{\,\wr_*\,}{{\mathbb Z}_2}$\\
\hline
6&$3K_2$&${{\mathbb Z}_2}\wr{S_3}$&${{\mathbb Z}_2}{\,\wr_*\,}{S_3}$\\
\hline
8&$2K_4$&$S_4\wr{{\mathbb Z}_2}$&$\mathcal Q_4{\,\wr_*\,}{{\mathbb Z}_2}$\\
\hline
8&$2C_4$&$({\mathbb Z}_2\wr{\mathbb Z}_2)\wr{{\mathbb Z}_2}$&$({\mathbb Z}_2{\,\wr_*\,}{\mathbb Z}_2){\,\wr_*\,}{{\mathbb Z}_2}$\\
\hline
8&$4K_2$&${{\mathbb Z}_2}\wr S_4$&${{\mathbb Z}_2}{\,\wr_*\,}\mathcal Q_4$\\
\hline
9&$3K_3$&${S_3}\wr{S_3}$&${S_3}{\,\wr_*\,}{S_3}$\\
\hline
10&$2C_5$&$D_5\wr{{\mathbb Z}_2}$&$D_5{\,\wr_*\,}{{\mathbb Z}_2}$\\
\hline
10&$2K_5$&$S_5\wr{{\mathbb Z}_2}$&$\mathcal Q_5{\,\wr_*\,}{{\mathbb Z}_2}$\\
\hline
10&$5K_2$&${{\mathbb Z}_2}\wr S_5$&${{\mathbb Z}_2}{\,\wr_*\,}\mathcal Q_5$\\
\hline
10&$C_{10}(4)$&${{\mathbb Z}_2}\wr D_5$&${{\mathbb Z}_2}{\,\wr_*\,} D_5$\\
\hline
\end{tabular}\end{center}

\medskip

\section{Classification table}

We are now in position of writing down a big table. We first recall the graph notations used in the paper.

\begin{definition}

We use the following notations.
\begin{enumerate}
\item Basic graphs:
- the complete graph having $n$ vertices is denoted $K_n$.

- the disjoint union of $n$ copies of $X$ is denoted $nX$.

- the prism having basis $X$ is denoted $\mathtt{Pr}(X)$.

\item Circulant graphs:

- the $n$-cycle is denoted $C_n$.

- the $2n$-cycle with diagonals is denoted $C_{2n}^+$.

- the $n$-cycle with chords of length $k$ is denoted $C_n(k)$.

\item Special graphs:

- the triangle times itself is denoted $\mathtt{Torus}$.

- the Petersen graph is denoted $\mathtt{Petersen}$.
\end{enumerate}
\end{definition}

As for quantum group notations, these have to be taken with care, because quantum groups do not really exist etc. Here they are.

\begin{definition}
We use the following notations.

- ${\mathbb Z}_n,D_n,S_n$ are the cyclic, dihedral and symmetric groups.

- $\mathcal Q_n$ is the quantum permutation group.

- $\times,\wr,{\,\wr_*\,}$ are the product, wreath product and free wreath product.
\end{definition}

The vertex-transitive graphs of order less than 11, modulo complementation, are given by the following table.

\vfill\eject

\begin{center}\begin{tabular}[t]{|l|l|l|l|}

\hline
Order & Graph & Classical group & Quantum group\\  
\hline \hline
2&$K_2$&$ {{\mathbb Z}_2}$&$ {{\mathbb Z}_2}$\\ 
\hline\hline
3&$K_3$&${{S_3}}$&${{S_3}}$\\ 
\hline\hline
4 & $2K_2$& ${\mathbb Z}_2\wr{\mathbb Z}_2$ & ${{\mathbb Z}_2}{\,\wr_*\,} {{\mathbb Z}_2}$ \\
\hline
4 & $K_4$ & $S_4$ & $\mathcal Q_4$   \\
\hline
\hline
5 & $C_5$ & $D_5$ & $D_5$ \\
\hline
5 & $K_5$ & $S_5$ & $\mathcal Q_5$\\
\hline \hline
6 & $C_6$ & $D_6$ & $D_6$ \\
\hline
6 & $2K_3$ & ${{S_3}}\wr{{{\mathbb Z}_2}}$ & ${{S_3}}{\,\wr_*\,}{{{\mathbb Z}_2}}$ \\
\hline
6 & $3K_2$ & ${{{\mathbb Z}_2}}\wr {{S_3}}$ & ${{{\mathbb Z}_2}}{\,\wr_*\,} {{S_3}}$ \\ 
\hline
6 & $K_6$ & $S_6$ & $\mathcal Q_6$ \\
\hline \hline
7 & $C_7$ & $D_7$ & $D_7$ \\
\hline
7 & $K_7$ & $S_7$ & $\mathcal Q_7$\\
\hline \hline
8 & $C_8$, $C_8^+$&  $D_8$ &  $D_8$\\
\hline
8 & $\texttt{Pr}(C_4)$ & 
$S_4 \times {{{\mathbb Z}_2}}$ & $\mathcal Q_4\times{{{\mathbb Z}_2}}$ \\
\hline
8 & $2K_4$ & $S_4\wr{{{\mathbb Z}_2}}$ &  $\mathcal Q_4{\,\wr_*\,}{{{\mathbb Z}_2}}$ \\
\hline
8 & $2C_4$& $({\mathbb Z}_2\wr{\mathbb Z}_2)\wr{{{\mathbb Z}_2}}$
& $({{{\mathbb Z}_2}}{\,\wr_*\,}{{{\mathbb Z}_2}}){\,\wr_*\,}{{{\mathbb Z}_2}}$ \\ 
\hline
8 & $4K_2$& ${{{\mathbb Z}_2}}\wr S_4$ & ${{{\mathbb Z}_2}}{\,\wr_*\,} \mathcal Q_4$ \\ 
\hline
8 & $K_8$ &  $S_8$ & $\mathcal Q_8$ \\
\hline \hline
9 & $C_9$, $C_9(3)$ & $D_9$ & $D_9$ \\
\hline
9 & $\texttt{Torus}$& ${{S_3}}\wr{{{\mathbb Z}_2}}$ & ${{S_3}}\wr{{{\mathbb Z}_2}}$ \\
\hline 
9 & $3K_3$ & ${{S_3}}\wr {{S_3}}$ & ${{S_3}}{\,\wr_*\,} {{S_3}}$ \\
\hline
9 & $K_9$ & $S_9$ & $\mathcal Q_9$ \\
\hline \hline
10 & $C_{10}$, $C_{10}(2)$, $C_{10}^+$, $\texttt{Pr}(C_5)$ & $D_{10}$ & $D_{10}$\\
\hline
10 & $\texttt{Petersen}$ & $S_5$ & $?$ \\
\hline
10 & $\texttt{Pr}(K_5)$ & $S_5 \times {{{\mathbb Z}_2}}$ &
$\mathcal Q_5\times{{{\mathbb Z}_2}}$ \\
\hline
10 & $C_{10}(4)$& ${{\mathbb Z}_2}\wr D_5$ & ${{\mathbb Z}_2}{\,\wr_*\,} D_5$\\
\hline
10 & $2C_5$ & $D_5\wr{{{\mathbb Z}_2}}$ & $D_5{\,\wr_*\,}{{{\mathbb Z}_2}}$ \\
\hline
10 & $2K_{5}$ & $S_5\wr{{{\mathbb Z}_2}}$ & $\mathcal Q_5{\,\wr_*\,}{{{\mathbb Z}_2}}$\\
\hline 
10 & $5K_2$ & ${{{\mathbb Z}_2}}\wr S_5$ & ${{{\mathbb Z}_2}}{\,\wr_*\,}\mathcal Q_5$ \\
\hline
10 & $K_{10}$ & $S_{10}$ &  $\mathcal Q_{10}$   \\
\hline \hline
11 & $C_{11}$, $C_{11}(2)$, $C_{11}(3)$& $D_{11}$ & $D_{11}$  \\
\hline
11 & $K_{11}$ & $S_{11}$ & $\mathcal Q_{11}$ \\
\hline
\end{tabular}\end{center}

\vfill\eject

Here the first three columns are well-known, and can be found in various books or websites. The last one collects results in this paper.

By using the equality $D_n={\mathbb Z}_n\rtimes {\mathbb Z}_2$, we reach the conclusion in the abstract: with one possible
exception, all quantum groups in the right column can be obtained from ${\mathbb Z}_n,S_n,\mathcal Q_n$ by using the operations $\times,\rtimes,\wr,{\,\wr_*\,}$.

The exceptional situation is that of the Petersen graph,
which might give a new quantum group.
We discuss it in the next section.

\section{The Petersen graph} 

The techniques of the previous sections do not apply to the Petersen graph,
which is not a circulant graph and cannot be written as a graph product.
Also we could not carry a  direct analysis similar to the one of the torus
because of the complexity of some computations.
The usual symmetry group is $S_5$, so
in view of the results in our classification table, we have at least two natural
candidates for the quantum symmetry group of the Petersen graph: $S_5$ and $\mathcal Q_5$. 

\begin{theorem}
The quantum automorphism group of the Petersen graph has an irreducible 
5-dimensional representation. In particular it is not 
isomorphic to the quantum symmetric group $\mathcal Q_5$.  
\end{theorem}

\begin{proof}
Let $G$ be the quantum automorphism group of the Petersen graph, denoted here $\texttt{P}$. We have an inclusion $S_5 \subset G$. It is well-known
that 
$${\rm Sp}(\mathtt{P}) = \{4, -2 , 1 \}$$
and that the corresponding eigenspaces have dimensions 
$1,4,5$. These eigenspaces furnish representations of
$G$ and of $S_5$. It is straightforward to compute the character
of the permutation representation of $S_5$ on $\mathbb C(\texttt{P})$, and then
using the character table of $S_5$
(see e.g. \cite{fh}), we see 
that $\mathbb C(\texttt{P})$ is the direct sum of 
$3$ irreducible representations of $S_5$. These have to be the previous
eigenspaces, and in particular the $5$-dimensional one is an irreducible
representation of $S_5$, and  of $G$.
On the other hand, it is known from \cite{ba0} 
that $\mathcal Q_5$ has no irreducible 5-dimensional representation. Thus the quantum groups $G$ and $\mathcal Q_5$ are not
isomorphic.
\end{proof}

The question now is: does the Petersen graph have quantum symmetry?
In other words, is $A(\texttt{P})$ commutative? 
The above result seems to indicate that if $A(\texttt{P})$ is not commutative,
we probably will have a new quantum permutation group.

\end{document}